\documentclass[11pt]{article}

\usepackage{amsmath,amssymb}
\usepackage[table]{xcolor}% http://ctan.org/pkg/xcolor

\usepackage{url}
\usepackage{graphicx}
\usepackage{subfigure}
\usepackage{color}
\usepackage{cite}
\usepackage{epsfig}
\usepackage{amssymb}
\usepackage{latexsym}
\newcommand\given[1][]{\:#1\vert\:}
\usepackage{mathtools}
\DeclarePairedDelimiter\ceil{\lceil}{\rceil}

\newcommand {\C} {{\rm I\kern-5.5pt C}}

\newcommand{\bP}[1]{{\mathbb{P}}\left[{#1}\right]}

\newcommand{\fsquare}{\vrule height6pt width7pt depth1pt}   % filled square
\newcommand{\myproof}{{\hfill \\ \bf Proof. \ }}           % Proof
\newcommand{\myendpf}{\hfill\fsquare \\[0.1in]}             % end of proof

\newtheorem{theorem}{Theorem}[section]

\newtheorem{proposition}[theorem]{Proposition}
\newtheorem{corollary}[theorem]{Corollary}

\newtheorem{fact}[theorem]{Fact}

\usepackage{cite}

\newcommand{\norm}[1]{\left\lVert#1\right\rVert}

\newcommand\blfootnote[1]{%
  \begingroup
  \renewcommand\thefootnote{}\footnote{#1}%
  \addtocounter{footnote}{-1}%
  \endgroup
}

\allowdisplaybreaks
\begin{document}

\title{Connectivity of Inhomogeneous Random K-out Graphs}

\author{
Rashad~Eletreby and
        Osman~Ya\u{g}an \\
{\tt reletreby@cmu.edu}, {\tt oyagan@ece.cmu.edu} \\
Department of Electrical and Computer Engineering and CyLab\\
Carnegie Mellon University.
}

% make the title area
\maketitle
\blfootnote{A preliminary version of some of the material was presented at the IEEE Conference on Decision and Control in 2018 \cite{eletrebycdc2018} and at the IEEE International Symposium on Information Theory in 2019 \cite{eletrebyISIT_2019}. This work has been supported in part by the National Science Foundation through grant CCF-1617934. 
}

\begin{abstract}
We propose inhomogeneous random K-out graphs $\mathbb{H}(n; \pmb{\mu}, \pmb{K}_n)$, where each of the $n$  nodes is assigned to one of $r$ classes independently with a probability distribution $\pmb{\mu} = \{\mu_1, \ldots, \mu_r\}$. In particular, each node is classified as class-$i$ with probability $\mu_i>0$, independently. Each class-$i$ node {\em selects}  $K_{i,n}$  distinct nodes uniformly at random from among all other nodes. A pair of nodes are adjacent in $\mathbb{H}(n; \pmb{\mu}, \pmb{K}_n)$ if at least one selects the other. Without loss of generality, we assume that $K_{1,n} \leq K_{2,n} \leq \ldots \leq K_{r,n}$. Earlier results on {\em homogeneous} random K-out graphs $\mathbb{H}(n; K_n)$, where all nodes select the same number $K$ of other nodes, reveal that $\mathbb{H}(n; K_n)$ is connected with high probability (whp) if $K_n \geq 2$ which implies that $\mathbb{H}(n; \pmb{\mu}, \pmb{K}_n)$ is connected whp if $K_{1,n} \geq 2$. In this paper, we investigate the connectivity of inhomogeneous random K-out graphs $\mathbb{H}(n; \pmb{\mu}, \pmb{K}_n)$ for the special case when $K_{1,n}=1$, i.e., when each class-$1$ node selects only one other node. We show that $\mathbb{H}\left(n;\pmb{\mu},\pmb{K}_n\right)$ is connected whp if $K_{r,n}$ is chosen such that $\lim_{n \to \infty} K_{r,n} = \infty$. However, any bounded choice of the sequence $K_{r,n}$ gives a positive probability of $\mathbb{H}\left(n;\pmb{\mu},\pmb{K}_n\right)$ being {\em not} connected. Simulation results are provided to validate our results in the finite node regime.
\end{abstract}


\section{Introduction}
\label{sec:Introduction}

The study of random graphs in their own right dates back to 1959 with the seminal work of Paul Erd\H{o}s, Alfred R\'enyi, and Edgar Gilbert. In particular, Erd\H{o}s and R\'enyi \cite{ER} introduced the random graph model $\mathbb{G}(n; M)$, representing a graph selected uniformly at random from the collection of all graphs with $n$ nodes and $M$ edges. In the same year, Gilbert \cite{gilbert1959} independently introduced the random graph model $\mathbb{G}(n; p)$, where each pair of vertices is connected (respectively, not connected) by an edge independently with probability $p$ (respectively, $1-p$). Since then, random graphs have received great attention in their own right and as a modeling framework for a wide class of real-world networks including social networks, communication networks, biological networks, among others \cite{boccaletti2006complex, goldenberg2010survey, newman2002random}.

Over the past several decades, a large body of research has developed several random graph models and investigated their structural properties \cite{Bollobas, JansonLuczakRucinski, frieze2016introduction}. For instance, {\em random geometric graphs} \cite{PenroseBook} $\mathbb{G}(n;r)$ were used to model the wireless connectivity of wireless ad-hoc networks \cite{Gupta99, haenggi2009stochastic}, providing guidelines on the critical radius needed to ensure network connectivity. In general, random geometric graphs are relevant to the modeling of networks that have a {\em spatial content} such as in wireless communications, epidemiology, and the internet \cite{PenroseBook}. Random geometric graphs are constructed as follows. Consider a set of vertices $\mathcal{V}=\{1,2,\ldots,n\}$, respectively located at random positions which are independent and uniformly distributed over a bounded region $\mathcal{D}$ of a Euclidean plane. Two vertices $i$ and $j$ located at $\pmb{x}_i$ and $\pmb{x}_j$, respectively, are connected by an edge if $\norm{\pmb{x}_i-\pmb{x}_j}<r$.

{\em Random key graphs} \cite{yagan2012zero, Blackburn_2009, DiPietroTissec} $\mathbb{G}(n;K,P)$ represent another class of random graphs that are used to model a wide range of applications and real-world networks such as common-interest friendship networks \cite{Yagan2017_SmallWorld} as well as secure connectivity of wireless sensor networks utilizing Eschenauer and Gligor random key predistribution scheme \cite{Gligor_2002, DiPietroTissec, yagan2012zero}\footnote{In Eschenauer and Gligor random key predistribution scheme \cite{Gligor_2002}, each sensor node is given (before deployment) $K$ cryptographic keys selected uniformly at random from a large key pool of size $P$. After deployment, two sensor nodes can communicate securely if they share a common key and are within wireless communication range.}. Random key graphs are constructed as follows. Each of the $n$ nodes is given $K$ objects selected uniformly at random (without replacement) from an object pool of size $P$. An edge exists between two vertices $u$ and $v$ only if they share an object. In the context of common-interest friendship networks, objects represent hobbies or interests, and two individuals are considered friends if they share a hobby or interest. In the context of secure connectivity of wireless sensor networks, objects represent cryptographic keys, and two sensor nodes can establish secure communication, only if they share a key.

Of particular interest to this paper is another commonly studied class of random graphs known as {\em random K-out graphs} $\mathbb{H}(n;K)$ \cite{Bollobas, Yagan2013Pairwise, FennerFrieze1982} that are constructed as follows. Each of the $n$ nodes selects $K$ other nodes uniformly at random from among all other others. An {\em undirected} edge is assigned between nodes $u$ and $v$ if $u$ selects $v$ or $v$ selects  $u$, or both; see \cite{Yagan2013Pairwise} for details. Random K-out graphs have recently received great interest for their role in modeling secure connectivity of wireless sensor networks utilizing Chan et al. \cite{Haowen_2003, Yagan2013Pairwise,yaugan2013scalability} random pairwise key predistribution scheme\footnote{In Chan et al. \cite{Haowen_2003} random pairwise key predistribution scheme, each of the $n$ sensor nodes is paired (offline) with $K$ distinct nodes which are randomly selected from among all other nodes. If nodes $i$ and $j$ were paired during the node-pairing stage, a unique (pairwise) key is generated and stored in the memory modules of each of the paired sensors together with both their IDs. After deployment, a secure link can be established between two communicating nodes if they have at least one pairwise key in common.}. More recently, a structure similar to random K-out graphs was suggested by Fanti et al. \cite[Algorithm~1]{FantiDandelion2018} to provide anonymity guarantees for transactions over cryptocurrency networks. The connectivity of random K-out graphs was studied in \cite{FennerFrieze1982, Yagan2013Pairwise}, where it was shown that
\begin{equation} 
\lim_{n \to \infty} \mathbb{P}\left[ \mathbb{H}(n;K) \text{ is connected}\right] =
\begin{cases}
0 & \mathrm{if} \quad K=1 \\
1 & \mathrm{if} \quad K\geq 2
\end{cases}
\label{eq:introEq}
\end{equation}
Hence, it is sufficient to set $K=2$ to have a connected network with high probability in the limit of large network size. In fact, it was shown in \cite{Yagan2013Pairwise} that the probability of $ \mathbb{H}(n;2)$ being connected exceeds 0.99 with as little as $n=50$ nodes. 

The aforementioned random graph models could all be described as {\em homogeneous} models, due to their uniform treatment of all vertices. In particular, homogeneous random graph models inherently assume that all vertices are {\em similar}, e.g., each vertex samples the same number $K$ of objects in random key graphs and each vertex selects the same number $K$ of other vertices in random K-out graphs. However, real-world complex networks are essentially composed of heterogeneous entities \cite{Barabasi_1999, boccaletti2006complex}, inducing the need for {\em inhomogeneous} variants of classical random graph models. Emerging wireless sensor networks represent a pronounced example of heterogeneous networks that consist of different nodes with different levels of resources (for communication, computation, storage, power, etc.) and possibly a varying level of security and connectivity requirements \cite{Du2007_applications, Lu2008_applications, Wu2007_applications, Yarvis_2005}. In fact, the literature on random graphs is already shifting towards inhomogeneous models initiated by the seminal work of Bollob\'as et al. on inhomogeneous Erd\H{o}s-Re\'nyi graph \cite{bollobas2007phase} (see also \cite{devroye2014connectivity}).

In this paper, we propose {\em inhomogeneous} random K-out graphs $\mathbb{H}(n; \pmb{\mu}, \pmb{K})$, where each of the $n$ nodes is assigned to one of $r$ classes according to a probability distribution $\pmb{\mu}=\{\mu_1, \ldots, \mu_r\}$ with $\mu_i>0$ for $i=1,\ldots,r$. A class-$i$ node selects $K_{i,n}$ nodes uniformly at random from among all other nodes. Two nodes $u$ and $v$ are connected by an edge if $u$ selects $v$ or $v$ selects $u$, or both. Without loss of generality, we assume that $K_{1,n} \leq K_{2,n} \leq \ldots \leq K_{r,n}$. We let $K_{\mathrm{avg},n}=\sum_{i=1}^r \mu_i K_{i,n}$ denote the expected number of selections. Inhomogeneous random K-out graphs generalize standard random K-out graphs to heterogeneous setting where different nodes make different number of selections depending on their corresponding classes. As a result, it might be expected that inhomogeneous K-out graphs would serve as a more natural model in many of the envisioned applications of K-out graphs including pairwise key predistribution in sensor networks and anonymous transactions in cryptocurrency networks.

By an easy monotonicity argument, we see from  (\ref{eq:introEq}) that $\mathbb{H}(n; \pmb{\mu}, \pmb{K}_n)$ is connected with high probability if $K_{1,n} \geq 2$. Of particular interest to our paper is the special case when $K_{1,n}=1$, i.e., when each of the nodes belonging to class-$1$ selects only one other node to be paired to. One could reasonably conjecture that setting $K_{\mathrm{avg},n} $ to any finite number larger than or equal to two would be sufficient to ensure the connectivity of $\mathbb{H}(n; \pmb{\mu}, \pmb{K}_n)$, in resemblance to (\ref{eq:introEq}). Our results reveal that such a conjecture does not hold and that the connectivity of $\mathbb{H}(n; \pmb{\mu}, \pmb{K}_n)$ under the special case when $K_{1,n}=1$ cannot be inferred from (\ref{eq:introEq}).

In this paper, we study the connectivity of  $\mathbb{H}(n; \pmb{\mu}, \pmb{K}_n)$ when $K_{1,n}=1$. More precisely, we seek conditions on $K_{2,n}, K_{3,n}, \ldots, K_{r,n}$ and $\pmb{\mu}$ such that the resulting graph is connected with high probability. Our main results (see Theorems \ref{thm:ZeroLaw+Connectivity} and \ref{thm:OneLaw+Connectivity}) show that $\mathbb{H}\left(n;\pmb{\mu},\pmb{K}_n\right)$ is connected with high probability if and only if $K_{r,n}= \omega(1)$. In other words, if $K_{r,n}$ grows unboundedly large as $n \to \infty$, then the probability 
that $\mathbb{H}(n; \pmb{\mu}, \pmb{K}_n)$ is connected approaches one in the same limit.
However, any bounded choice of $K_{r,n}$ gives a positive probability of $\mathbb{H}\left(n;\pmb{\mu},\pmb{K}_n\right)$ being {\em not} connected in the limit of large $n$. Comparing our results with (\ref{eq:introEq}) sheds the light on a striking difference between inhomogeneous random K-out graphs $\mathbb{H}\left(n;\pmb{\mu},\pmb{K}_n\right)$ and their homogeneous counterpart $\mathbb{H}\left(n;K\right)$. In particular, the flexibility of organizing the nodes into several classes with different characteristics (with $K_{1,n}=1$) comes at the expense of requiring $\lim_{n \to \infty} K_{r,n} = \infty$, in contrast to the homogeneous case where having $K=2$ was sufficient to ensure connectivity.

Throughout the paper, all statements involving limits, including asymptotic equivalences, are understood with $n$ going to infinity. The cardinality of any discrete set $S$ is denoted by $|S|$. The random variables (rvs) under consideration are all defined on the same probability triple $(\Omega, {\cal F}, \mathbb{P})$. Probabilistic statements are made with respect to this probability measure $\mathbb{P}$, and we denote the corresponding expectation operator by $\mathbb{E}$. We say that an event holds with high probability (whp) if it holds with probability $1$ as $n \rightarrow \infty$.  In comparing the asymptotic behaviors of the sequences $\{a_n\},\{b_n\}$, we use $a_n = o(b_n)$,  $a_n=\omega(b_n)$, and $a_n = O(b_n)$ with their meaning in the standard Landau notation. %Namely, we write  $a_n = o(b_n)$ (respectively, $a_n = \omega(b_n)$) as a shorthand for the relation $\lim_{n \to \infty} \frac{a_n}{b_n}=0$ (respectively, $\lim_{n \to \infty} \frac{a_n}{b_n}=\infty$), whereas $a_n = O(b_n)$ means that there exists $c>0$ such that $a_n \leq c b_n$ for all $n$ sufficiently large. 
We write $\mathbb{N}_0$ to denote the set of natural numbers excluding zero.

\section{Inhomogeneous random K-out graphs}
\label{sec:Model}

The inhomogeneous random K-out graph, denoted $\mathbb{H}\left(n; \pmb{\mu}, \pmb{K}_n\right)$, is constructed on the vertex set $\mathcal{V}=\{1,2,\ldots,n\}$ as follows. First, each node is assigned a class $i \in \{1,\ldots, r\}$ independently according to a probability distribution
$\pmb{\mu} = \{\mu_1, \ldots, \mu_r \}$; i.e.,  $\mu_i$ denotes the probability that a node is class-$i$ and we have  $\sum_{i=1}^r \mu_i=1$. We assume $\mu_i>0$ for all $i=1,2,\ldots,r$ and that $r$ is a fixed integer that does not scale with %the number of nodes 
$n$. Each class-$i$ node selects $K_{i,n}$ distinct nodes uniformly at random from $\mathcal{V}\setminus \{v\}$ and an undirected edge is assigned between a pair of nodes if at least one selects the other.  Formally, each node $v$ is associated  (independently from others) with a subset $\Gamma_{n,v}(\pmb{\mu},\pmb{K}_n)$ (whose size depends on the class of node $v$) of nodes selected uniformly at {\em random} from $\mathcal{V}\setminus \{ v \}$.
%Each of the nodes in $\Gamma_{n,v}(\pmb{\mu},\pmb{K}_n)$ is
%said to be {\em paired} to node $i$.
Specifically, for any  $A \subseteq \mathcal{V}\setminus \{ v \}$, we have
\begin{equation}
\bP{ \Gamma_{n,v}(\pmb{\mu},\pmb{K}_n)= A \given[\big] t_v=i } = \left \{
\begin{array}{ll}
{{n-1}\choose{K_i}}^{-1} & \mbox{if $|A|=K_i$} \\
              &                   \\
0             & \mbox{otherwise}
\end{array}
\right .
\label{eq:main_eqn_for_gamma}
\end{equation}
where $t_v$ denotes the class of node $v$. Then, vertices $u$ and $v$ are said to be adjacent in $\mathbb{H}\left(n; \pmb{\mu}, \pmb{K}_n\right)$, written $u \sim v$, if 
at least one selects the other; i.e., if
%$u$ selects $v$ or $v$ selects $u$, or both.  Namely
\begin{equation}
u \sim v
\quad \mbox{iff} \quad
u \in \Gamma_{n,v} (\pmb{\mu},\pmb{K}_n) ~{\huge \vee}~ v \in \Gamma_{n,u}(\pmb{\mu},\pmb{K}_n).
\label{eq:Adjacency}
\end{equation}

When $r=1$, all vertices belong to the same class and thus select the same number (say, $K$) of other nodes, leading to the {\em homogeneous} random K-out graph  $\mathbb{H}(n;K)$ \cite{Bollobas,FennerFrieze1982, Yagan2013Pairwise}.

Throughout, we set
\begin{equation}
K_{\mathrm{avg},n} = \sum_{i=1}^r \mu_i K_{i,n}
\label{eq:kAvg_def}
\end{equation}

For any distinct nodes $u, v \in \mathcal{V}$, we have
\begin{align}
\mathbb{P}\left[ u \sim v \right] &= 1 - \mathbb{P}\left[ u \not\in \Gamma_{n,v} (\pmb{\mu},\pmb{K}_n) \cap  v \not\in \Gamma_{n,u} (\pmb{\mu},\pmb{K}_n)\right]= 1- \left(	\sum_{i=1}^r \mu_i \frac{\binom{n-2}{K_i}}{\binom{n-1}{K_i}} 	\right)^2=1 - \left(1- \frac{K_{\mathrm{avg},n}}{n-1}\right)^2 \nonumber
\end{align}

\section{Main results}
\label{sec:Results}
We refer to any mapping $\pmb{K}: \mathbb{N}_0
\rightarrow \mathbb{N}^r_0$ as a {\em scaling} provided it satisfies the condition.
\begin{equation} 
K_{1,n} \leq K_{2,n} \leq \ldots \leq K_{r,n} < n, \quad n=2,3, \ldots . 
\label{eq:scalingCond}
\end{equation}

Our main technical results, given next, characterize the connectivity of inhomogeneous random K-out graphs. Throughout, it will be convenient to use the notation
\begin{equation} \nonumber
P(n;\pmb{\mu},\pmb{K}_n):= \mathbb{P} \left[ \mathbb{H}(n;\pmb{\mu},\pmb{K}_n) \text{ is connected} \right]
\end{equation}
and
\begin{equation} 
C(\pmb{\mu},\pmb{K}_n) = \frac{1}{1+\frac{2}{\mu_1^2}e^{2 K_{\mathrm{avg},n}}}
\label{eq:defining_Ck}
\end{equation}
and
\begin{align}
\Psi(n, \pmb{\mu}, \pmb{K}_n) &=  \max \Big\{ \exp \left( -2\left( 1-\tilde{\mu} \right) \left(\frac{K_{r,n}-1}{4}  - \frac{(0.5)^{K_{r,n}-1}}{\tilde{\mu}} \right) \right) , \nonumber \\
&  \mathrm{exp} \left( -\left( 1-\tilde{\mu} \right) \frac{n}{2} \left( 1-e^{-1}  - \frac{(0.5)^{K_{r,n}-1}}{\tilde{\mu}} \right)  \right) \Big\} 
\label{eq:defining_Psi}
\end{align}
with $0 < \mu_1 <1$, $K_{\mathrm{avg},n}$ as defined in (\ref{eq:kAvg_def}), and $\tilde{\mu} = \sum_{i=1}^{r-1} \mu_i$.

The following result establishes an upper bound on the probability of connectivity of the inhomogeneous random K-out graphs when the sequence $K_{r,n}$ is bounded, i.e., $K_{r,n}=O(1)$

\begin{theorem}
{\sl  Consider a scaling $\pmb{K}: \mathbb{N}_0 \rightarrow \mathbb{N}^r_0$ and a probability distribution $\pmb{\mu}=\{ \mu_1,\mu_2, \ldots, \mu_r \}$ with $\mu_i>0$.  %such that a node is labeled as class-$i$ with probability $\mu_i$. 
If $K_{r,n} = O(1)$, then
\begin{equation}
\limsup_{n \to \infty} P(n;\pmb{\mu},\pmb{K}_n) <1
\label{eq:zeroLawStatement}
\end{equation}
More precisely, we have
\begin{equation}
P(n;\pmb{\mu},\pmb{K}_n) \leq 1-C(\pmb{\mu},\pmb{K}_n)  + o(1)
\label{eq:zeroLawPrecise}
\end{equation}
}
\label{thm:ZeroLaw+Connectivity}
\end{theorem}

The following result establishes a one-law for connectivity for the inhomogeneous random K-out graph. %s when $K_{r,n} = \omega(1)$.
\begin{theorem}
{\sl Consider a scaling $\pmb{K}: \mathbb{N}_0 \rightarrow \mathbb{N}^r_0$ and a probability distribution $\pmb{\mu}=\{\mu_1,\mu_2, \ldots, \mu_r\}$ with $\mu_i>0$. %such that a node is labeled as class-$i$ with probability $\mu_i$. 
If $K_{r,n} = \omega(1)$, then
\begin{equation}\nonumber
\lim_{n \rightarrow \infty }P(n;\pmb{\mu},\pmb{K}_n) = 1
\end{equation}
More precisely, we have
\begin{equation}
P(n;\pmb{\mu},\pmb{K}_n) \geq 1- \frac{\tilde{\mu}^2 }{1-\tilde{\mu}} \Psi(n, \pmb{\mu}, \pmb{K}_n)
\label{eq:oneLawPrecise}
\end{equation}
for all $K_{r,n}$ sufficiently large such that
$ K_{r,n} \geq \ceil*{4 \left( \frac{(0.5)^{K_{r,n}-1}}{\tilde{\mu}}\right)+1}$.
}
\label{thm:OneLaw+Connectivity}
\end{theorem}

Theorems \ref{thm:ZeroLaw+Connectivity} and \ref{thm:OneLaw+Connectivity} state that $\mathbb{H}\left(n;\pmb{\mu},\pmb{K}_n\right)$ is connected with high probability if $K_{r,n}$ is chosen such that $K_{r,n} = \omega(1)$. On the other hand, if $K_{r,n} = O(1)$, then the probability of connectivity of $\mathbb{H}\left(n;\pmb{\mu},\pmb{K}_n\right)$ is strictly less than one in the limit of large network size. In other words, any {\em bounded} choice for $K_{r,n}$ gives rise to a positive probability of $\mathbb{H}\left(n;\pmb{\mu},\pmb{K}_n\right)$ being {\em not} connected. Observe that (\ref{eq:zeroLawStatement}) follows from (\ref{eq:zeroLawPrecise}) by virtue of the fact that $K_{\mathrm{avg},n} = O(1)$ when $K_{r,n} = O(1)$.

\subsection{Discussion}

Connectivity results in the literature of random graphs are usually presented in the form of zero-one laws, where the probability of connectivity (in the limit as $n \to \infty$) exhibits a sharp transition between two different regimes. In the first (respectively, second) regime, the probability tends to zero (respectively, one) as $n$ tends to infinity. One example of such results is given by (\ref{eq:introEq}) where the probability that $\mathbb{H}(n; K)$ is connected tends to zero when $K=1$ and tends to one when $K \geq 2$. Other examples include the connectivity results on random key graphs \cite{yagan2012zero}, Erd\H{o}s-R\'enyi graphs \cite{ER}, etc. Indeed, Theorem~\ref{thm:ZeroLaw+Connectivity} states that the probability of connectivity is strictly less than one whenever $K_{r,n}=O(1)$ but it does not specify whether or not a zero-law exists in this case. In other words, Theorem~\ref{thm:ZeroLaw+Connectivity} does not reveal whether or not $\lim_{n \to \infty} P(n;\pmb{\mu},\pmb{K}_n) =0$ when $K_{r,n}=O(1)$. Such a zero-law, if exists, would complement the one-law given by Theorem~\ref{thm:OneLaw+Connectivity}.

A careful look at (\ref{eq:oneLawPrecise}) reveals that $P(n;\pmb{\mu},\pmb{K}_n)$ exhibits a lower bound that could either be trivial (negative) or non-trivial (positive). As a result, under the conditions that force the bound to be non-trivial, the probability of connectivity is strictly larger than zero, hence a zero-law does not exist in this case. In what follows, we let $K^\star(\tilde{\mu})$ denote the smallest value of $K_{r,n}$ for which 
\begin{equation}\nonumber
\frac{\tilde{\mu}^2 }{1-\tilde{\mu}} \Psi(n, \pmb{\mu}, \pmb{K}_n) < 1.
\end{equation}
We present a result that utilizes (\ref{eq:oneLawPrecise}) to show that under some conditions on $\tilde{\mu}$ and $K_{r,n}$, the probability of connectivity of $\mathbb{H}\left(n;\pmb{\mu},\pmb{K}_n\right)$ is strictly larger than zero, hence, a zero-law does not hold. 

{
\corollary
\label{cor:cor1}
Consider a scaling $\pmb{K}: \mathbb{N}_0 \rightarrow \mathbb{N}^r_0$ and a probability distribution $\pmb{\mu}=\{\mu_1,\mu_2, \ldots, \mu_r\}$ with $\mu_i>0$. For any $\tilde{\mu}$, there exists $K^\star(\tilde{\mu})$ such that 
\begin{equation}\nonumber
P(n;\pmb{\mu},\pmb{K}_n)>0
\end{equation}
whenever $K_{r,n} \geq K^\star(\tilde{\mu})$.
}

In Table~\ref{table:t1}, we provide the values of $K^\star(\tilde{\mu})$ corresponding to some values of $\tilde{\mu}$. Note that whether or not a zero-law holds for the case when 
$2 \leq K_{r,n} < K^\star(\tilde{\mu})$ cannot by established through (\ref{eq:oneLawPrecise}) and is beyond the scope of this paper.

\begin{table}
\begin{center}
 \begin{tabular}{||c | c||| c | c ||} 
 \hline
 \cellcolor{gray!15} $\tilde{\mu}$ &  \cellcolor{gray!15} $K^\star(\tilde{\mu})$ & \cellcolor{gray!15} $\tilde{\mu}$ & \cellcolor{gray!15} $K^\star(\tilde{\mu})$ \\ [0.5ex] 
 \hline\hline
 0.1 & 5 & 0.6 & 3  \\ 
 \hline
 0.2 & 4  & 0.7 & 5  \\ 
 \hline
 0.3 & 4 &  0.8 & 13  \\ 
 \hline
 0.4 & 4  & 0.9 & 43\\
 \hline
 0.5 & 3 & 0.95 & 117 \\ [1ex] 
 \hline
\end{tabular}
\caption{The values of $K^\star(\tilde{\mu})$ corresponding to different values for $\tilde{\mu}$. When $K_{r,n} \geq K^\star(\tilde{\mu})$, the probability of connectivity of $\mathbb{H}\left(n;\pmb{\mu},\pmb{K}_n\right)$ is strictly larger than zero by virtue of (\ref{eq:oneLawPrecise}), hence a zero-law does not hold in this case.}
\label{table:t1}
\end{center}
\end{table}

Theorems~\ref{thm:ZeroLaw+Connectivity} and \ref{thm:OneLaw+Connectivity} reveal a striking difference between inhomogeneous random K-out graphs and their homogeneous counterpart. In particular, it was shown in \cite{FennerFrieze1982, Yagan2013Pairwise} that
\begin{equation} \nonumber
\lim_{n \to \infty} \mathbb{P}\left[ \mathbb{H}(n;K) \text{ is connected}\right] =
\begin{cases}
0 & \mathrm{if} \quad K=1 \\
1 & \mathrm{if} \quad K\geq 2
\end{cases}
\end{equation}
Hence, it is sufficient to set $K=2$ to have a connected network with high probability in the limit of large network size. When the network size $n$ is fixed, Ya\u{g}an and Makowski \cite{Yagan2013Pairwise} showed that
\begin{equation}\nonumber
\mathbb{P}\left[ \mathbb{H}(n;2) \text{ is connected}\right] \geq 1-\frac{155}{n^3}, \qquad n\geq16
\end{equation} 
indicating that the probability of connectivity exceeds $0.99$ for as little as $n=50$ nodes (with $K=2$). As a result, random K-out graphs $\mathbb{H}(n;K)$ can be connected with orders of magnitude fewer links, in total, as compared to most other random graph models such as Erd\H{o}s-R\'enyi graphs \cite{ER}, random key graphs \cite{yagan2012zero}, and inhomogeneous random key graphs \cite{Yagan/Inhomogeneous}, where the mean degree (respectively, the {\em minimum} mean degree in inhomogeneous random key graphs) has to be on the order of $\log n$ to ensure connectivity. In contrast, the mean degree of $\mathbb{H}\left(n;K\right)$ is of order $2 K$, i.e., a mean degree of $4$ is sufficient to ensure connectivity of $\mathbb{H}(n;K)$.

In contrast, {\em inhomogeneous} random K-out graphs (with $K_{1,n}=1$) require $K_{r,n}$ to grow unboundedly large as $n \to \infty$ so that the probability of connectivity approaches one in the same limit. In other words, the flexibility of arranging nodes into classes comes at the expense of {\em sparsity}. In particular, the mean degree of $\mathbb{H}\left(n; \pmb{\mu}, \pmb{K}\right)$ has to grow unboundedly large as $n \to \infty$ to ensure the connectivity of the graph. Fortunately, Theorem~\ref{thm:OneLaw+Connectivity} does not specify a particular growth rate function for the sequence $K_{r,n}$, other than $K_{r,n}=\omega(1)$. Hence, one can set $K_{r,n} = \log \log \ldots \log n$ to meet the requirements of Theorem~\ref{thm:OneLaw+Connectivity}. As a result, inhomogeneous random K-out graphs $\mathbb{H}(n;K)$ can be connected with orders of magnitude fewer links, in total, as compared to most other random graph models as mentioned above.

%A special case of the inhomogeneous random K-out graph was investigated in \cite{eletrebycdc2018} where the number of classes was limited to two, i.e., $r=2$. Therein, each node is classified as class-$1$ with probability $\mu$ or class-$2$ with probability $1-\mu$. Each of class-$1$ (respectively, class-$2$) nodes selects $1$ (respectively, $K_n$) nodes chosen uniformly at random from the rest of the nodes. In contrast to \cite{eletrebycdc2018}, we consider a {\em generalization} of inhomogeneous random K-out graphs where the number of classes is not necessarily limited to two. In particular, we consider the general case where each node belongs to one of $r$ classes, and each class-$i$ node selects $K_{i,n}$ other nodes chosen uniformly at random (with $K_{1,n}=1$). Our proof techniques 
%Even though the proof of Theorem~\ref{thm:ZeroLaw+Connectivity} in \cite{eletrebycdc2018} is limited to the case where $r=2$, it is superfluous and requires $K_{n}$ to admit a particular value $K$ for all $n=2,3,\ldots$. In contrast, our proof handles the general case of $r$ different classes and does not require such stringent conditions on $K_{i,n}$ for $i=2,\ldots,r$. Namely, in establishing Theorem~\ref{thm:ZeroLaw+Connectivity} we do not require any of the sequences $K_{i,n}$ to admit a particular value for all $n=2,3,\ldots$, or to have a particular limit. In addition to generalizing the proof and relaxing the conditions of \cite{eletrebycdc2018}, we also use different proof techniques that lead to a more concise and sharper proof.

\subsection{Numerical Study}
\label{subsec:numerical}

\begin{figure}[!t]
\vspace{-5mm}
\hspace{-.3cm}
\centering
\includegraphics[scale=0.55]{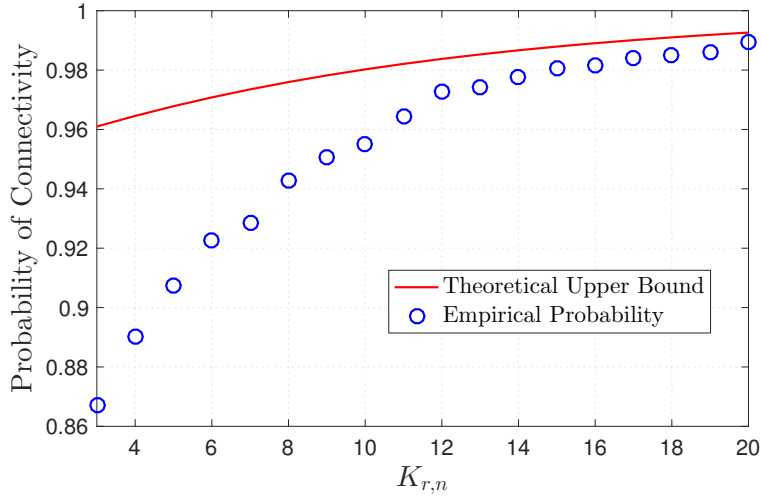} 
\vspace{-2mm}
\caption{\sl The empirical probability $P(n;\pmb{\mu},\pmb{K}_n)$ with $\pmb{\mu}=\left\{0.9, 0.06, 0.04\right\}$ and $\pmb{K}=(1, 2, K_3)$ as a function of $K_{3}$ for $n=1000$ along with the theoretical upper bound given by Theorem~\ref{thm:ZeroLaw+Connectivity}. Empirical probabilities approach the upper bound as $K_{3}$ increases. Empirical probabilities were obtained by averaging over $10^5$ independent experiments for each data point.}
\label{fig:prob}
\end{figure}

The objective of this subsection is to validate the upper bound given by Theorem~\ref{thm:ZeroLaw+Connectivity} in the finite-node regime using computer simulations. In Figure~\ref{fig:prob}, we consider an inhomogeneous random K-out graph with three classes. Namely, we set $\pmb{\mu} = \left\{0.9, 0.06, 0.04 \right\}$ and $\pmb{K}=(1, 2, K_3)$, i.e., each node is classified as class-$1$ with probability $0.9$, class-$2$ with probability $0.06$, and class-$3$ with probability $0.04$. Nodes belonging to class-$1$ (respectively, class-$2$) select only one (respectively, two) other node(s) to be paired to. We vary $K_3$ from $3$ to $20$ and observe how the empirical probability of connectivity varies in accordance. In particular, for each value of $K_3$, we run $10^5$ independent experiments for each data point and count the number of times (out of $10^5$) when the resulting graph is connected. Dividing this number by $10^5$ gives the {\em empirical} probability of connectivity.

Note that as $K_3$ varies, $K_{\mathrm{avg}}$ varies as well according to (\ref{eq:kAvg_def}). We can then use (\ref{eq:defining_Ck}) to plot the theoretical upper bound given by $1-C(\pmb{\mu}, \pmb{K})$. The results given in Figure \ref{fig:prob} confirm the validity of Theorem~\ref{thm:ZeroLaw+Connectivity} but also reveals its shortcomings. Observe that the bound appears to be loose for small values of $K_{3}$, yet it becomes tighter as $K_3$ increases. The reasoning behind this phenomenon would become apparent in Section~\ref{sec:proofTh1} as we outline our approach in establishing Theorem~\ref{thm:ZeroLaw+Connectivity}. At a high level, our approach is based on bounding the probability of connectivity by the probability of {\em not} observing isolated components of size two, i.e., components formed by two class-$1$ nodes $u$ and $v$ such that $u$ has selected $v$, $v$ has selected $u$, and none of the other $n-2$ nodes has either selected $u$ or $v$. When $K_3$ is large, the probability of observing isolated components of sizes larger than two (i.e., three, four, etc.) will be small. Hence, the probability of connectivity in this regime would be tightly bounded by the probability of {\em not} observing isolated components of size two. However, in the regime where $K_3$ is small, isolated components of sizes other than two are more likely to be formed, as compared to the case when $K_3$ is large (see Figure~\ref{fig:realization}). Since our approach does not consider such components, our bound becomes slightly loose in this regime.

\section{Preliminaries}
\label{section:prelim}
Throughout, we will make use of the following results.

{
\fact [{\cite[Fact~2]{Jun/K-Connectivity}}]
\label{fact1xynew}
For $0 \leq x < 1$, and $y=0,1,2,\ldots$, we have
\begin{equation} \nonumber
1-xy\leq \left(1-x\right)^y \leq 1-xy+\frac{1}{2}x^2y^2
\end{equation}
}

{
\fact [{\cite[Fact~4]{Jun/K-Connectivity}}]
\label{PncKn}
Let integers $x$ and $y$ be both positive
functions of $n$, where $y \geq 2x$. For $z = 0, 1, \ldots, x$, we have
\begin{align}
\frac{\binom{y- z}{x} } {\binom{y}{x}} & \geq 1 - \frac{ z x }{y - z
}, \label{fone}
\end{align}
and
\begin{align}
 \frac{\binom{y- z}{x} } {\binom{y}{x}}
 & = 1 - \frac{xz}{y}
 \pm O\bigg(\frac{x^4}{y^2}\bigg). \label{fone2}
\end{align}
}

{
{
\fact 
For $r=1,\ldots,\lfloor \frac{n}{2} \rfloor $ and $n=1,2,\ldots$, we have
\begin{equation}
\binom{n}{r}  \leq \left( \frac{n}{r}\right)^r \left(\frac{n}{n-r}\right)^{n-r}
\label{eq:Bound_on_NchooseR}
\end{equation}
}

\myproof
The following bound, established in \cite{robbins1955remark}, is valid for all $x=1, 2, \ldots$ 
\begin{equation}
\sqrt{2 \pi} x^{x+0.5} e^{-x} e^{\frac{1}{12x+1}} < x! < \sqrt{2 \pi} x^{x+0.5} e^{-x} e^{\frac{1}{12x}}.
\label{eq:robbinsBound}
\end{equation}
Observe that
\begin{equation}\nonumber
\sqrt{2 \pi} e^{\frac{1}{12x}} \leq e
\end{equation}
for all $x\geq2$.
and
\begin{equation}\nonumber
 e^{\frac{1}{12x+1}} \geq 1
\end{equation}
 Hence, (\ref{eq:robbinsBound}) can be written as
\begin{equation}
\sqrt{2 \pi} x^{x+0.5} e^{-x} < x! < e x^{x+0.5} e^{-x}
\label{eq:robbinsBound_v2}
\end{equation}
Using (\ref{eq:robbinsBound_v2}), we get
\begin{align}
\binom{n}{r} &=\dfrac{n!}{r! (n-r)!} \nonumber \\
&< \dfrac{en^{n+0.5}e^{-n}}{\sqrt{2 \pi} r^{r+0.5} e^{-r} \sqrt{2 \pi} (n-r)^{n-r+0.5}e^{-(n-r)}} \nonumber \\
&= \dfrac{e}{2 \pi } \dfrac{1}{ \sqrt{r} \sqrt{1-\frac{r}{n}}} \dfrac{1}{\left( \frac{r}{n}\right)^r \left( 1-\frac{r}{n}\right)^{n-r}} \nonumber \\
& \leq \dfrac{e}{2 \pi  \sqrt{0.5}} \dfrac{1}{\left( \frac{r}{n}\right)^r \left( 1-\frac{r}{n}\right)^{n-r}} \nonumber \\
& \leq \left( \frac{n}{r}\right)^r \left(\frac{n}{n-r}\right)^{n-r}
\end{align}
as we use the crude bounds $r \geq 1$ and $r \leq n/2$.
\myendpf
}

For $0 \leq K \leq x \leq y$, we have
\begin{equation}
\dfrac{\binom{x}{K}}{\binom{y}{K}} = \prod_{\ell=0}^{K-1} \left(\frac{x-\ell}{y-\ell} \right) \leq \left( \frac{x}{y} \right)^K
\label{eq:boundOnChooseK}
\end{equation}
since $\frac{x-\ell}{y-\ell}$ decreases as $\ell$ increases from $\ell=0$ to $\ell=K-1$.

Moreover, we have
\begin{equation}
1 \pm x \leq e^{\pm x}, \quad 0 \leq x \leq 1
\label{eq:expBound} \\
\end{equation}
and
\begin{equation}
1 - e^{- x} \geq \frac{x}{2}, \quad 0 \leq x \leq 1
\label{eq:expBound2} \\
\end{equation}
%
%Finally, we find it useful to write
%\begin{equation}
%1-x=e^{-x-\Psi(x)}, \quad 0 \leq x <1
%\label{eq:isolated_log_decomp}
%\end{equation}
%where  $\Psi(x)=\int_{0}^{x} \frac{t}{1-t} \ \text{dt}$.
%From L'H\^{o}pital's Rule, we have
%\begin{equation}
%\lim_{x\to 0}  \frac{\Psi(x)}{x^2}=\frac{-x-\log (1-x)}{x^2}=\frac{1}{2}.
%\label{eq:isolated_hopital}
%\end{equation}

Throughout, we set 
\begin{equation} 
\binom{x}{y}=0,
\label{eq:prelim:zeroBinom}
\end{equation}
whenever $x<y$.

\section{A proof of Theorem~\ref{thm:ZeroLaw+Connectivity}}
\label{sec:proofTh1}

In what follows, we establish (\ref{eq:zeroLawPrecise}) whenever $K_{n,r} = O(1)$. In particular, with each class-$1$ node selecting only {\em one} other node, we will show that whenever each class-$r$ node gets paired to a {\em bounded} number of nodes, there will be a positive probability that the graph is {\em not} connected. Note that if the sequence $K_{r,n}$ is bounded, then so are the sequences $K_{i,n}$  for $i=2,\ldots,r-1$ by virtue of (\ref{eq:scalingCond}). Put differently
\begin{equation}\nonumber
K_{r,n} = O(1) \Rightarrow K_{i,n} = O(1), \quad i=2,\ldots,r-1
\end{equation}

\begin{figure}[!t]
\vspace{-5mm}
\hspace{-.3cm}
\centering
\includegraphics[scale=0.045]{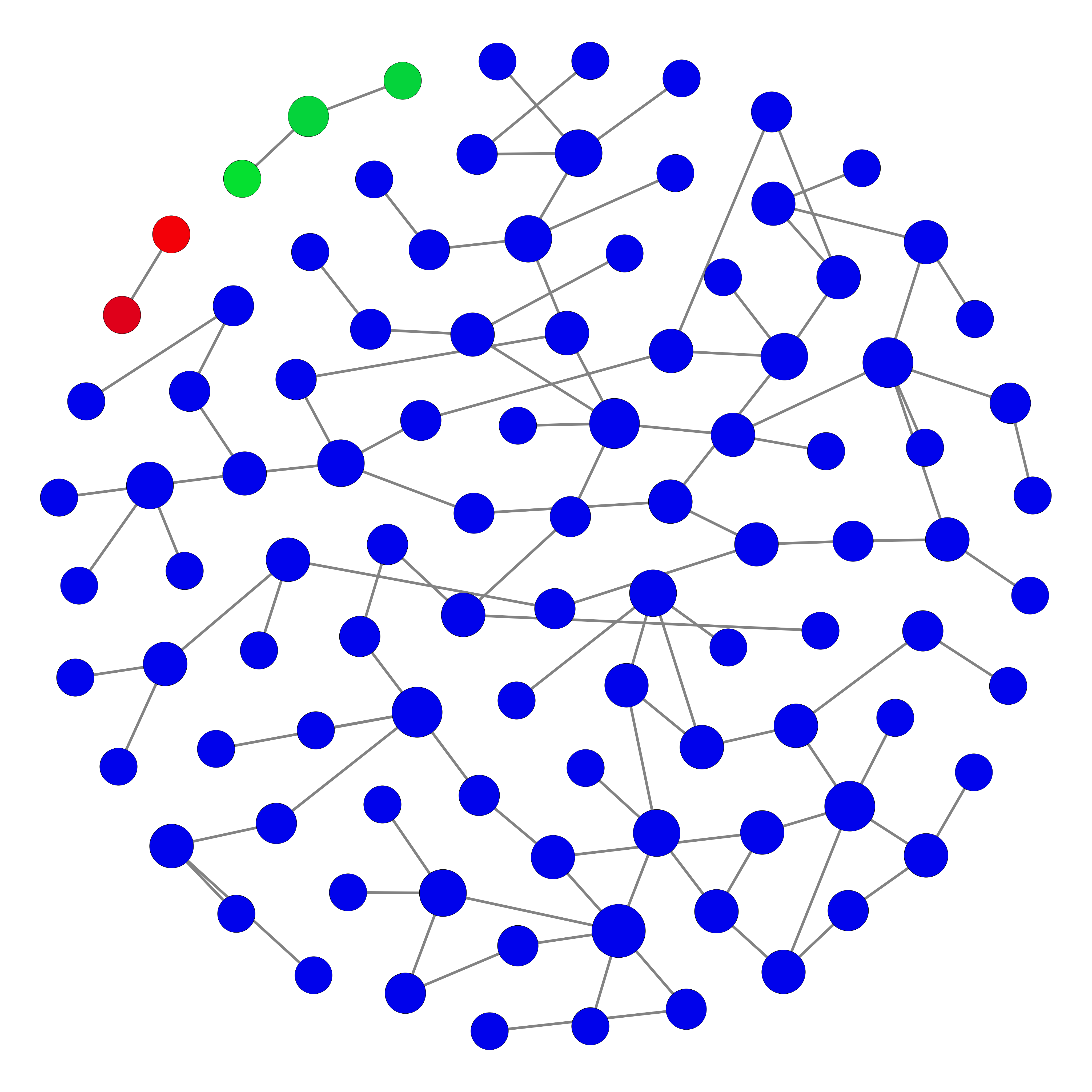} 
\vspace{-2mm}
\caption{\sl A realization of the inhomogeneous random K-out graph $\mathbb{H}(n;\pmb{\mu},\pmb{K}_n)$ with $r=3$, $\pmb{K}=\left(1, 2, 3\right)$. The graph is {\em not} connected as it contains two {\em isolated} components, highlighted in red and green, respectively. The first isolated component consists of two nodes, while the second isolated component consists of three nodes. We set $n=100$ and $\pmb{\mu}=\left\{0.9, 0.05, 0.05\right\}$. The size of each node corresponds to its degree. \vspace{-3mm}}
\label{fig:realization}
\end{figure}

Observe that when a positive fraction of the nodes, each, gets paired with only one node, the graph may contain isolated components consisting of two class-$1$ nodes, say $i$ and $j$, that were paired with each other, i.e., $\Gamma_{n,i}(\pmb{\mu},\pmb{K}_n)=\{j\}$, $\Gamma_{n,j}(\pmb{\mu},\pmb{K}_n)=\{i\}$, and $\Gamma_{n,\ell}(\pmb{\mu},\pmb{K}_n) \subseteq \mathcal{N} \setminus \{i,j,\ell\}$ for all $\ell \in \mathcal{N}\setminus\{i,j\}$. Indeed, these isolated components render the graph disconnected. A graphical illustration is given in Figure~\ref{fig:realization}. Our approach in establishing Theorem~\ref{thm:ZeroLaw+Connectivity} relies on the method of second moment applied to a variable that counts the number of isolated components that contain two vertices of class-$1$. 

Recall that $t_i$ denotes the class of node $i$. Let $U_{ij}(n;\pmb{\mu},\pmb{K}_n)$ denote the event that nodes $i$ and $j$ are both class-$1$ {\em and} are forming an isolated component, i.e.,
\begin{align}
&U_{ij}(n;\pmb{\mu},\pmb{K}_n) \label{eq:defining_U} \\
& = \left( \bigcap_{\ell \in \mathcal{N}\setminus \{i,j\}} \left[ \Gamma_{n,\ell}(\pmb{\mu},\pmb{K}_n) \subseteq \mathcal{N} \setminus \{i,j,\ell\} \right] \right) \cap 
\left[\Gamma_{n,i}(\pmb{\mu},\pmb{K}_n)=\{j\} \right]  \cap \left[ \Gamma_{n,j}(\pmb{\mu},\pmb{K}_n)=\{i\} \right] \cap \left[t_1=1 \right] \cap \left[t_2=1 \right] \nonumber
\end{align}
Next, let 
\begin{align}
\chi_{ij}(n;\pmb{\mu},\pmb{K}_n) = \pmb{1} \left[ U_{ij}(n;\pmb{\mu},\pmb{K}_n)\right] \nonumber
\end{align}
and
\begin{equation}\nonumber
Y(n;\pmb{\mu},\pmb{K}_n) = \sum_{1 \leq i < j \leq n} \chi_{ij}(n;\pmb{\mu},\pmb{K}_n)
\end{equation}

Clearly, $Y(n;\pmb{\mu},\pmb{K}_n)$ gives the number of isolated components in $\mathbb{H}(n;\pmb{\mu},\pmb{K}_n)$ that contain two vertices of class-$1$. We will show that when $K_{r,n} = O(1)$, we have
\begin{equation}\nonumber
\mathbb{P} \left[ Y(n;\pmb{\mu},\pmb{K}_n)=0 \right] \leq 1-C(\pmb{\mu},\pmb{K}_n) + o(1)
\end{equation}
Recall that if $\mathbb{H}(n;\pmb{\mu},\pmb{K}_n)$ is connected, then it does not contain any isolated component. In particular,  $\mathbb{H}(n;\pmb{\mu},\pmb{K}_n)$ would consist of a single component of size $n$. However, the absence of isolated components of size two does not necessarily mean that $\mathbb{H}(n;\pmb{\mu},\pmb{K}_n)$  is connected, as it may contain isolated components of other sizes (see Figure~\ref{fig:realization}). It follows that,
\begin{align} 
P(n;\pmb{\mu},\pmb{K}_n) &\leq \mathbb{P} \left[ Y(n;\pmb{\mu},\pmb{K}_n)=0 \right] \nonumber 
\end{align}
Hence, establishing (\ref{eq:zeroLawPrecise}) is equivalent to establishing 
\begin{equation} 
\mathbb{P} \left[ Y(n;\pmb{\mu},\pmb{K}_n)=0 \right] \leq 1-C(\pmb{\mu},\pmb{K}_n) +o(1)
\label{eq:ZeroLaw_proof_equivalent}
\end{equation}
where $C(\pmb{\mu},\pmb{K}_n)$ is given by (\ref{eq:defining_Ck}).

By applying the method of second moments
\cite[Remark 3.1, p. 55]{JansonLuczakRucinski}  on $Y(n;\pmb{\mu},\pmb{K}_n)$, we get
\begin{equation} 
\mathbb{P}[Y(n;\pmb{\mu},\pmb{K}_n)=0] \leq 1-\frac{\left(\mathbb{E}[Y(n;\pmb{\mu},\pmb{K}_n)]\right)^2}{\mathbb{E}[Y^2(n;\pmb{\mu},\pmb{K}_n)]}  
\label{eq:upperBoundSecondMoment}
\end{equation}
where
\begin{align} 
\mathbb{E}[Y(n;\pmb{\mu},\pmb{K}_n)] &= \sum_{1 \leq i < j \leq n} \mathbb{E} \left[\chi_{ij}(n;\pmb{\mu},\pmb{K}_n)\right] =\binom{n}{2} \mathbb{E}[\chi_{12}(n;\pmb{\mu},\pmb{K}_n)] 
\label{eq:isolated_ZeroLaw_first_part}
\end{align}
and
\begin{align}
\mathbb{E}[Y^2(n;\pmb{\mu},\pmb{K}_n)] &= \mathbb{E} \left[\sum_{1 \leq i < j \leq n} \sum_{1 \leq \ell < m \leq n} \chi_{ij}(n;\pmb{\mu},\pmb{K}_n) \chi_{\ell m}(n;\pmb{\mu},\pmb{K}_n)  \right]\nonumber \\
&= \binom{n}{2} \mathbb{E} \left[ \chi_{12}(n;\pmb{\mu},\pmb{K}_n) \right] + 2\binom{n}{2} \binom{n-2}{1} \mathbb{E} \left[ \chi_{12}(n;\pmb{\mu},\pmb{K}_n) \chi_{13}(n;\pmb{\mu},\pmb{K}_n) \right] \nonumber \\
& \quad + \binom{n}{2} \binom{n-2}{2} \mathbb{E} \left[ \chi_{12}(n;\pmb{\mu},\pmb{K}_n) \chi_{34}(n;\pmb{\mu},\pmb{K}_n) \right] \nonumber 
\end{align}
by exchangeability and the binary nature of the random variables $\{\chi_{ij}(n;\pmb{\mu},\pmb{K}_n) \}_{1 \leq i <j \leq n}$. Observe that
\begin{equation}\nonumber
\mathbb{E} \left[ \chi_{12}(n;\pmb{\mu},\pmb{K}_n) \chi_{13}(n;\pmb{\mu},\pmb{K}_n) \right]  = 0,
\end{equation} 
since $\left[ U_{12}(n;\pmb{\mu},\pmb{K}_n) \cap U_{13}(n;\pmb{\mu},\pmb{K}_n) \right]=\emptyset$ by definition. Hence,
\begin{align}
\mathbb{E}[Y^2(n;\pmb{\mu},\pmb{K}_n)] = \binom{n}{2} \mathbb{E} \left[ \chi_{12}(n;\pmb{\mu},\pmb{K}_n) \right] + \binom{n}{2} \binom{n-2}{2} \mathbb{E} \left[ \chi_{12}(n;\pmb{\mu},\pmb{K}_n) \chi_{34}(n;\pmb{\mu},\pmb{K}_n) \right]  \label{eq:isolated_ZeroLaw_second_part} 
\end{align}

Using (\ref{eq:isolated_ZeroLaw_first_part}) and  (\ref{eq:isolated_ZeroLaw_second_part}), we get
\begin{align}
&\frac{\mathbb{E}[Y^2(n;\pmb{\mu},\pmb{K}_n)]}{\left(\mathbb{E}[Y(n;\pmb{\mu},\pmb{K}_n)]\right)^2} = \frac{1}{\binom{n}{2} \mathbb{E}[\chi_{1,2}(n;\pmb{\mu},\pmb{K}_n)]}   + \frac{\binom{n}{2} \binom{n-2}{2} \mathbb{E}[\chi_{1,2}(n;\pmb{\mu},\pmb{K}_n) \chi_{3,4}(n;\pmb{\mu},\pmb{K}_n)]}{\left( \binom{n}{2} \mathbb{E}[\chi_{1,2}(n;\pmb{\mu},\pmb{K}_n)]\right)^2}
\label{eq:isolated_ZeroLaw_second_part_detailed}
\end{align}

The next two results will help establish (\ref{eq:ZeroLaw_proof_equivalent}). 

{
\proposition
\sl
\label{prop:firstMoment}
Consider a scaling $\pmb{K}: \mathbb{N}_0 \rightarrow \mathbb{N}^r_0$ and a probability distribution $\pmb{\mu}=\{ \mu_1,\mu_2, \ldots, \mu_r \}$ with $\mu_i>0$.  %such that a node is labeled as class-$i$ with probability $\mu_i$. 
If $K_{r,n} = O(1)$, then
\begin{align}
\hspace{-2mm} \binom{n}{2} \mathbb{E}\left[ \chi_{12} (n;\pmb{\mu},\pmb{K}_n)\right] =   \left(1+o(1) \right) \frac{\mu_1^2}{2} \mathrm{exp} \left(-2 K_{\mathrm{avg},n}\right)   \label{eq:zeroLaw_firstProp}
\end{align}
}

\myproof
Note that under $U_{12}(n;\pmb{\mu},\pmb{K}_n)$, we have 
\begin{equation}\nonumber
\Gamma_{n,1}(\pmb{\mu},\pmb{K}_n)=\{2\} \quad \mathrm{and} \quad \Gamma_{n,2}(\pmb{\mu},\pmb{K}_n)=\{1\}
\end{equation}
Moreover, we have
\begin{equation}\nonumber
\Gamma_{n,i}(\pmb{\mu},\pmb{K}_n) \subseteq \mathcal{N} \setminus \{1,2,i\}, \quad i=3,4,\ldots,n
\end{equation}
Recall that each of the other $n-2$ nodes is class-$i$ with probability $\mu_i$ and that the random variables $\Gamma_{n,1}(\pmb{\mu},\pmb{K}_n),\Gamma_{n,2}(\pmb{\mu},\pmb{K}_n),\ldots,\Gamma_{n,n}(\pmb{\mu},\pmb{K}_n)$ are mutually independent. Hence, we have
\begin{align}
\mathbb{E}\left[ \chi_{12} (n;\pmb{\mu},\pmb{K}_n) \right] &= \mathbb{P}\left[ U_{12}(n;\pmb{\mu},\pmb{K}_n) \right] = \mu_1^2 \left( \frac{1}{n-1}	\right)^2 \left( \sum_{i=1}^r \mu_i \frac{\binom{n-3}{K_{i,n}}}{\binom{n-1}{K_{i,n}}}   \right)^{n-2}  \nonumber
\end{align}
Then, we have
\begin{align}
 \binom{n}{2} \mathbb{E}\left[ \chi_{12} (n;\pmb{\mu},\pmb{K}_n)\right] &= \frac{\mu_1^2}{2} \left( \frac{n}{n-1}	\right)  \left( \sum_{i=1}^r \mu_i \frac{\binom{n-3}{K_{i,n}}}{\binom{n-1}{K_{i,n}}}   \right)^{n-2}  \nonumber
\\ 
&= \frac{\mu_1^2}{2} \left( \frac{n}{n-1}	\right) \cdot  \left( \sum_{i=1}^r \mu_i \left( \frac{\left( n-1-K_{i,n} \right)}{\left(n-1 \right)} \frac{\left( n-2-K_{i,n} \right)}{\left(n-2 \right)}  \right) \right)^{n-2}  \nonumber
\\
&= \frac{\mu_1^2}{2} \left( \frac{n}{n-1}	\right)  \left( \sum_{i=1}^r \mu_i \left( 1-\frac{K_{i,n}}{n-1}  \right) \left( 1-\frac{K_{i,n}}{n-2} \right) \right)^{n-2}  \nonumber 
\\ 
&= \frac{\mu_1^2}{2} \left( \frac{n}{n-1}	\right)  \cdot  \left( 1 - \left(\sum_{i=1}^r \mu_i  \frac{2K_{i,n}(n-1.5)}{(n-1)(n-2)} \right)+\left(\sum_{i=1}^r \mu_i  \frac{K^2_{i,n}}{(n-1)(n-2)}\right)  \right)^{n-2}  \nonumber 
\\
& =  \frac{\mu_1^2}{2} \left( \frac{n}{n-1}	\right)  \cdot \exp\left(-2 \left(\frac{n-1.5}{n-1}\right) \sum_{i=1}^r \mu_i  K_{i,n} + \frac{1}{n-1}  \sum_{i=1}^r \mu_i  K^2_{i,n} \right)
\nonumber
\\
&= \left(1 + o(1) \right) \frac{\mu_1^2}{2} e^{-2 K_{\mathrm{avg},n}} \nonumber
\end{align}
where the last equality follows since $K_{r,n} = O(1)$.
\myendpf

{
\proposition
\sl
Consider a scaling $\pmb{K}: \mathbb{N}_0 \rightarrow \mathbb{N}^r_0$ and a probability distribution $\pmb{\mu}=\{ \mu_1,\mu_2, \ldots, \mu_r \}$ with $\mu_i>0$.  %such that a node is labeled as class-$i$ with probability $\mu_i$. 
If $K_{r,n} = O(1)$, then
\begin{align}
&  \frac{  \mathbb{E}\left[ \chi_{12} (n;\pmb{\mu},\pmb{K}_n) \chi_{34} (n;\pmb{\mu},\pmb{K}_n)\right] }{\left( \mathbb{E}\left[ \chi_{12} (n;\pmb{\mu},\pmb{K}_n)\right] \right)^2 }  = 1 +o(1).
\label{eq:zeroLaw_secondProp}
\end{align}
}

\myproof
Note that an immediate consequence of Fact~\ref{PncKn} is that
\[
\frac{\binom{n}{2} \binom{n-2}{2}}{\binom{n}{2} ^2} = 1+o(1)
\]

Observe that under $\left[ U_{12}(n;\pmb{\mu},\pmb{K}_n) \cap U_{34}(n;\pmb{\mu},\pmb{K}_n) \right]$, we have 
\begin{align}
&\Gamma_{n,1}(\pmb{\mu},\pmb{K}_n)=\{2\} \quad \mathrm{and} \quad \Gamma_{n,2}(\pmb{\mu},\pmb{K}_n)=\{1\} \nonumber \\
&\Gamma_{n,3}(\pmb{\mu},\pmb{K}_n)=\{4\} \quad \mathrm{and} \quad \Gamma_{n,4}(\pmb{\mu},\pmb{K}_n)=\{3\} \nonumber 
\end{align}
Moreover, we have
\begin{equation}\nonumber
\Gamma_{n,i}(\pmb{\mu},\pmb{K}_n) \subseteq \mathcal{N} \setminus \{1,2,3,4,i\}, \quad i=5,6,\ldots,n
\end{equation}
Recall that each of the other $n-4$ nodes is class-$i$ with probability $\mu_i$ and that the random variables $\Gamma_{n,1}(\pmb{\mu},\pmb{K}_n),\Gamma_{n,2}(\pmb{\mu},\pmb{K}_n),\ldots,\Gamma_{n,n}(\pmb{\mu},\pmb{K}_n)$ are mutually independent. Hence, we have
\begin{align}
\mathbb{E}\left[ \chi_{12} (n;\pmb{\mu},\pmb{K}_n) \chi_{34} (n;\pmb{\mu},\pmb{K}_n)\right] = \mathbb{P} \left[U_{12}(n;\pmb{\mu},\pmb{K}_n) \cap U_{34}(n;\pmb{\mu},\pmb{K}_n)  \right] =\mu_1^4 \left(\frac{1}{n-1} \right)^4 \left( \sum_{i=1}^r \mu_i \frac{\binom{n-5}{K_{i,n}}}{\binom{n-1}{K_{i,n}}} \right)^{n-4} \nonumber 
\end{align}
Invoking Fact~\ref{PncKn}, we get
\begin{align}\nonumber
\frac{\binom{n}{2} \binom{n-2}{2} \mathbb{E}\left[ \chi_{12} (n;\pmb{\mu},\pmb{K}_n) \chi_{34} (n;\pmb{\mu},\pmb{K}_n)\right]}{\left(\binom{n}{2} \mathbb{E}\left[ \chi_{12} (n;\pmb{\mu},\pmb{K}_n)\right] \right)^2} &=  (1+o(1)) \frac{\left( \sum_{i=1}^r \mu_i \frac{\binom{n-5}{K_{i,n}}}{\binom{n-1}{K_{i,n}}} \right)^{n-4}}{\left( \sum_{i=1}^r \mu_i \frac{\binom{n-3}{K_{i,n}}}{\binom{n-1}{K_{i,n}}} \right)^{2n-4}} \nonumber \\
&=  (1+o(1)) \cdot \frac{\left( \sum_{i=1}^r \mu_i \left( 1 - \frac{4K_{i,n}}{n-1} \pm O \left( \frac{K^4_{i,n}}{n^2} \right) \right) \right)^{n-4}}{\left( \sum_{i=1}^r \mu_i \left( 1 - \frac{2K_{i,n}}{n-1} \pm O \left( \frac{K^4_{i,n}}{n^2} \right) \right) \right)^{2n-4}} \nonumber
\nonumber \\ 
& =  (1+o(1)) \cdot \frac{\left(1 - \frac{4 \sum_{i=1}^r \mu_i K_{i,n}}{n-1} \pm O \left( \frac{1}{n^2} \right) \right)^{n-4}}{\left(1 - \frac{2 \sum_{i=1}^r \mu_i K_{i,n}}{n-1} \pm O \left( \frac{1}{n^2} \right) \right)^{2n-4}} \nonumber
\nonumber \\ 
& =  (1+o(1)) \cdot \left( \frac{1 - \frac{4 K_{\mathrm{avg}, n}}{n-1} \pm O \left( \frac{1}{n^2} \right)}{\left(1 - \frac{2 K_{\mathrm{avg},n}}{n-1} \pm O \left( \frac{1}{n^2} \right)\right)^2} \right)^n 
\nonumber \\ 
& =  (1+o(1)) \cdot \left( \frac{1 - \frac{4 K_{\mathrm{avg},n}}{n-1} \pm O \left( \frac{1}{n^2} \right)}{1 - \frac{4 K_{\mathrm{avg},n}}{n-1} \pm O \left( \frac{1}{n^2} \right)} \right)^n 
\nonumber \\ 
& =  (1+o(1)) \cdot  \left( 1  \pm O\left(\frac{1}{n^2}\right)   \right)^n = 1+o(1).
\nonumber
\end{align}
\myendpf

The main result (\ref{eq:zeroLawPrecise}) now follows by virtue of (\ref{eq:ZeroLaw_proof_equivalent}) and (\ref{eq:upperBoundSecondMoment}) as we combine (\ref{eq:isolated_ZeroLaw_second_part_detailed}), (\ref{eq:zeroLaw_firstProp}), and (\ref{eq:zeroLaw_secondProp}). Observe that (\ref{eq:zeroLawStatement}) follows from (\ref{eq:zeroLawPrecise}) by virtue of the fact that $K_{\mathrm{avg},n} = O(1)$ when $K_{r,n} = O(1)$.

\section{A proof of Theorem~\ref{thm:OneLaw+Connectivity}}
\label{sec:proofTh2}
In what follows, we establish that
\begin{equation}
\lim_{n \rightarrow\infty} P(n;\pmb{\mu},\pmb{K}_n)  =1
\label{eq:OneLaw_proof_expression}
\end{equation}
whenever $K_{r,n} = \omega(1)$.

Observe that for any non-empty subset $S$ of nodes, i.e., $S \subseteq {\cal N}$,
we say that $S$ is {\em isolated} in
$\mathbb{H} (n;\pmb{\mu},\pmb{K}_n)$ if there are no edges in $\mathbb{H} (n;\pmb{\mu},\pmb{K}_n)$ between the nodes in $S$ and the nodes in the complement
$S^c = {\cal N} - S$. This is characterized by the event
$B_n (\pmb{\mu},\pmb{K}_n; S)$ given by
\begin{align}
B_n (\pmb{\mu},\pmb{K}_n; S) =
\bigcap_{i \in S} \bigcap_{j \in S^c}
\left (
\left [ i \not \in \Gamma_{n,j}(\pmb{\mu},\pmb{K}_n) \right ] 
\cap 
\left [ j \notin \Gamma_{n,i}(\pmb{\mu},\pmb{K}_n) \right ]
\right ). \nonumber
\end{align}

Note that if $\mathbb{H} (n;\pmb{\mu},\pmb{K}_n)$ is {\em not} connected, then there must exist a
non-empty subset $S$ of nodes which is isolated. Recall that each node in $\mathbb{H} (n;\pmb{\mu},\pmb{K}_n)$ is class-$i$ with probability $\mu_i$ and that $K_{1,n}=1$. Thus, we may observe isolated sets in $\mathbb{H} (n;\pmb{\mu},\pmb{K}_n)$ of cardinality\footnote{Note that if vertices $S$ form an isolated set then so do vertices $\mathcal{N}-S$, hence the sum need to be taken only until $\lfloor \frac{n}{2} \rfloor $.} $\ell=2,3,\ldots, \lfloor \frac{n}{2} \rfloor $. Thus, with $D_n (\pmb{\mu},\pmb{K}_n)$ denoting the event that $\mathbb{H}(n;\pmb{\mu},\pmb{K}_n)$ is connected, we have the inclusion
\begin{equation}
D_n (\pmb{\mu},\pmb{K}_n)^c
\subseteq
\cup_{S \in \mathcal{P}_n: ~2\leq  |S| \leq \lfloor \frac{n}{2} \rfloor} ~ B_n (\pmb{\mu},\pmb{K}_n; S)
\label{eq:BasicIdea}
\end{equation}
where $\mathcal{P}_n$ stands for the collection of all non-empty
subsets of ${\cal N}$. 
A standard union bound argument immediately gives
\begin{align}
\bP{ D_n(\pmb{\mu},\pmb{K}_n)^c }
\leq 
\sum_{ S \in \mathcal{P}_n: 2 \leq |S| \leq \lfloor \frac{n}{2} \rfloor }
\bP{ B_n (\pmb{\mu},\pmb{K}_n; S) } =
\sum_{\ell=2}^{ \lfloor \frac{n}{2} \rfloor }
\left ( \sum_{S \in \mathcal{P}_{n,\ell} } \bP{ B_n (\pmb{\mu},\pmb{K}_n; S)} \right )
\label{eq:BasicIdea+UnionBound}
\end{align}
where $\mathcal{P}_{n,\ell} $ denotes the collection of all subsets
of ${\cal N}$ with exactly $\ell$ elements.

For each $\ell=1, \ldots , n$, we simplify the notation by writing
$B_{n,\ell} (\pmb{\mu},\pmb{K}_n)= B_n (\pmb{\mu},\pmb{K}_n ; \{ 1, \ldots , \ell \} )$. Under the enforced
assumptions, exchangeability implies
\[
\bP{ B_n (\pmb{\mu},\pmb{K}_n ; S) } = \bP{ B_{n,\ell} (\pmb{\mu},\pmb{K}_n) }, \quad S \in
\mathcal{P}_{n,\ell}
\]
and the expression
\begin{equation}
\sum_{S \in \mathcal{P}_{n,\ell} } \bP{ B_n (\pmb{\mu},\pmb{K}_n ; S) } 
= {n \choose \ell} ~ \bP{ B_{n,\ell} (\pmb{\mu},\pmb{K}_n)} 
\label{eq:ForEach=r}
\end{equation}
follows since $|\mathcal{P}_{n,\ell} | = {n \choose \ell}$. Substituting
into (\ref{eq:BasicIdea+UnionBound}) we obtain the bounds
\begin{align}
\bP{ D_n (\pmb{\mu},\pmb{K}_n)^c } \leq \sum_{\ell=2}^{ \lfloor \frac{n}{2} \rfloor }
{n \choose \ell} ~ \bP{ B_{n,\ell}(\pmb{\mu},\pmb{K}_n) } .
\label{eq:BasicIdea+UnionBound2}
\end{align}

For each $\ell=2, \ldots , \lfloor \frac{n}{2} \rfloor$, it is easy to check that
\begin{align}
\mathbb{P}\left[B_{n,\ell}(\pmb{\mu},\pmb{K}_n) \right] =
\left(\sum_{i=1}^r \mu_i \frac{\binom{\ell-1}{K_{i,n}} }{ \binom{n-1}{K_{i,n}}} \right)^{\ell} \left(\sum_{i=1}^r \mu_i \frac{\binom{n-\ell-1}{K_{i,n}}}{\binom{n-1}{K_{i,n}} } \right)^{n-\ell}
\label{eq:B_n_r_K}
\end{align}

To see why this last relation holds, recall  that for nodes $\{1,\ldots, \ell\}$ 
to be isolated in $\mathbb{H}(n;\pmb{\mu},\pmb{K}_n)$, we need that 
(i) none of the sets $\Gamma_{n,1}(\pmb{\mu},\pmb{K}_n), \ldots, \Gamma_{n,\ell}(\pmb{\mu},\pmb{K}_n)$ 
contains an element from the set  $\{\ell+1, \ldots, n \}$; 
and (ii) none of the sets $\Gamma_{n,\ell+1}(\pmb{\mu},\pmb{K}_n), \ldots, \Gamma_{n,n}(\pmb{\mu},\pmb{K}_n)$  
contains an element from $\{ 1, \ldots, \ell\}$. 
More precisely, we must have
\[
\Gamma_{n, i} (\pmb{\mu},\pmb{K}_n) \subseteq \{1, \ldots, \ell\} \setminus \{i\}, 
\quad i=1, \ldots, \ell
\]
and
\[
\Gamma_{n, j}(\pmb{\mu},\pmb{K}_n)
\subseteq \{\ell+1, \ldots, n\} \setminus \{j\}, 
\quad j=\ell+1, \ldots, n.
\]
 Hence, the validity of (\ref{eq:B_n_r_K}) is now immediate from 
(\ref{eq:main_eqn_for_gamma}) and the mutual 
independence of the rvs $\Gamma_{n, 1}(\pmb{\mu},\pmb{K}_n), \ldots, \Gamma_{n, n} (\pmb{\mu},\pmb{K}_n)$.

We now establish that under the enforced assumptions of Theorem~\ref{thm:OneLaw+Connectivity}, we have
\begin{equation}\nonumber
\lim_{n \to \infty} \sum_{\ell=2}^{ \lfloor \frac{n}{2} \rfloor }
{n \choose \ell} ~ \bP{ B_{n,\ell}(\pmb{\mu},\pmb{K}_n) } = 0
\end{equation}
which in turn establishes Theorem~\ref{thm:OneLaw+Connectivity} by virtue of (\ref{eq:BasicIdea+UnionBound2}).

Note that the quantities
\begin{equation}\nonumber
\frac{\binom{\ell-1}{K_{i,n}}}{\binom{n-1}{K_{i,n}}} \quad \mathrm{and} \quad \frac{\binom{n-\ell-1}{K_{i,n}}}{\binom{n-1}{K_{i,n}}}
\end{equation}
are monotonically decreasing in $K_{i,n}$. We use (\ref{eq:boundOnChooseK}) and (\ref{eq:B_n_r_K}) to get
\begin{align}
\mathbb{P}\left[ B_{n,\ell}(\pmb{\mu},\pmb{K}_n) \right] &=
\left(\sum_{i=1}^{r-1} \mu_i \frac{\binom{\ell-1}{K_{i,n}} }{ \binom{n-1}{K_{i,n}}} + \mu_r \frac{\binom{\ell-1}{K_{r,n}} }{ \binom{n-1}{K_{r,n}}} \right)^{\ell} \cdot \left(\sum_{i=1}^{r-1} \mu_i \frac{\binom{n-\ell-1}{K_{i,n}}}{\binom{n-1}{K_{i,n}} } + \mu_r \frac{\binom{n-\ell-1}{K_{r,n}}}{\binom{n-1}{K_{r,n}} } \right)^{n-\ell} \nonumber \\
& \leq
\left(\frac{\binom{\ell-1}{K_{1,n}} }{ \binom{n-1}{K_{1,n}}} \left( \sum_{i=1}^{r-1} \mu_i  \right)+ \mu_r \frac{\binom{\ell-1}{K_{r,n}} }{ \binom{n-1}{K_{r,n}}} \right)^{\ell} \cdot \left( \frac{\binom{n-\ell-1}{K_{1,n}}}{\binom{n-1}{K_{1,n}} } \left(\sum_{i=1}^{r-1} \mu_i\right)+ \mu_r \frac{\binom{n-\ell-1}{K_{r,n}}}{\binom{n-1}{K_{r,n}} } \right)^{n-\ell} \nonumber \\
&\leq \left( \tilde{\mu} \left( \frac{\ell-1}{n-1}\right) + \left(1-\tilde{\mu}\right) \left( \frac{\ell-1}{n-1}\right)^{K_{r,n}} \right)^\ell  \cdot  \left( \tilde{\mu} \left( \frac{n-\ell-1}{n-1}\right) + \left(1-\tilde{\mu}\right) \left( \frac{n-\ell-1}{n-1}\right)^{K_{r,n}} \right)^{n-\ell} \label{eq:convertedBound}
\end{align}
where $\tilde{\mu} = \sum_{i=1}^{r-1} \mu_i$ and $1-\tilde{\mu} = \mu_r$. 

Observe that the bound appearing in (\ref{eq:convertedBound}) resembles the case where each node belongs to one of two classes. Namely, a node could either be class-$1$ (with probability $\tilde{\mu}$) or class $r$ (with probability $1-\tilde{\mu}$). We further use (\ref{eq:expBound}) to get
\begin{align}
\mathbb{P}\left[ B_{n,\ell}(\pmb{\mu},\pmb{K}_n ) \right] &\leq \left( \tilde{\mu} \left( \frac{\ell}{n} \right) + \left( 1-\tilde{\mu}  \right) \left( \frac{\ell}{n}\right)^{K_{r,n}} \right)^\ell \cdot \left(\tilde{\mu}  \left(1-\frac{\ell}{n} \right)+\left( 1-\tilde{\mu}  \right) \left(1-\frac{\ell}{n} \right)^{K_{r,n}} \right)^{n-\ell} \nonumber \\
&=\tilde{\mu}^\ell \left( \frac{\ell}{n} \right)^\ell  \left( 1 + \frac{1-\tilde{\mu} }{\tilde{\mu} } \left( \frac{\ell}{n}\right)^{K_{r,n}-1} \right)^\ell
\left( 1-\frac{\ell}{n} \right)^{n-\ell} \cdot \left( 1-\left( 1-\tilde{\mu}  \right) \left(1 - \left(1-\frac{\ell}{n} \right)^{K_{r,n}-1} \right) \right)^{n-\ell} \nonumber \\
&\leq \tilde{\mu}^\ell \left( \frac{\ell}{n} \right)^\ell  \left( 1-\frac{\ell}{n} \right)^{n-\ell} \left( 1 + \frac{1-\tilde{\mu} }{\tilde{\mu} } \left( \frac{\ell}{n}\right)^{K_{r,n}-1} \right)^\ell
  \cdot \left( 1-\left( 1-\tilde{\mu}  \right) \left(1 - e^{-\ell \left(\frac{K_{r,n}-1}{n}\right)} \right) \right)^{n-\ell} \nonumber \\
&\leq \tilde{\mu}^\ell \left( \frac{\ell}{n} \right)^\ell  \left( 1-\frac{\ell}{n} \right)^{n-\ell}  \mathrm{exp}\Bigg( \frac{1-\tilde{\mu} }{\tilde{\mu} }\ell \left( \frac{\ell}{n}\right)^{K_{r,n}-1}  -  \left( 1-\tilde{\mu}  \right) (n-\ell) \left( 1-e^{-\ell \left( \frac{K_{r,n}-1}{n}\right)}\right) \Bigg)  \label{eq:firstRangeBNR} 
\end{align}
Combining (\ref{eq:Bound_on_NchooseR}) with (\ref{eq:firstRangeBNR}), we conclude that
\begin{align}
 \mathbb{P} \left[D_n (\pmb{\mu},\pmb{K}_n)^c \right]
 \leq
 \sum_{\ell=2}^{ \lfloor \frac{n}{2} \rfloor }
 \binom{n}{\ell} \mathbb{P}\left[ B_{n,\ell}(\pmb{\mu},\pmb{K}_n ) \right]  \leq    \sum_{\ell=2}^{ \lfloor \frac{n}{2} \rfloor }  \tilde{\mu}^\ell A_{n,\ell}  
 \label{eq:last_bound_osy}
\end{align}
where we define
\begin{align}\label{eq:defn_A_nr}
 A_{n,\ell} & :=\mathrm{exp}\Bigg( \frac{1-\tilde{\mu}}{\tilde{\mu}}\ell \left( \frac{\ell}{n}\right)^{K_{r,n}-1}  -  \left( 1-\tilde{\mu} \right) (n-\ell) \left( 1-e^{-\ell \left( \frac{K_{r,n}-1}{n}\right)}\right) \Bigg)  
 \end{align}
with $2 \leq \ell \leq n/2$.

Next, our goal is to derive an upper bound on $A_{n,\ell}$ that is valid for all 
 $n$ sufficiently large and $\ell=2,\ldots, \lfloor \frac{n}{2} \rfloor$, and show that this bound tends to zero as 
$n$ gets large. 
Fix $n=2, 3,$ sufficiently large. For each $\ell=2, \ldots, \lfloor \frac{n}{2} \rfloor$,
either one of the following should hold
 \begin{equation}\nonumber
 \frac{\ell(K_{r,n}-1)}{n} \leq 1 \quad \mathrm{and} \quad \frac{\ell(K_{r,n}-1)}{n} > 1.
 \end{equation} 
 If it holds that $ \frac{\ell(K_{r,n}-1)}{n} \leq 1$, then we use 
 (\ref{eq:expBound2}) to get $1-e^{-\ell \left( \frac{K_{r,n}-1}{n}\right)} \geq  \frac{\ell (K_{r,n}-1)}{2n}$.
 Using this in (\ref{eq:defn_A_nr}) 
  yields
\begin{align}\nonumber
A_{n,\ell} &\leq \mathrm{exp} \left(\frac{1-\tilde{\mu}}{\tilde{\mu}} \ell \left(\frac{\ell}{n} \right)^{K_{r,n}-1} -\left( 1-\tilde{\mu} \right) (n-\ell) \frac{\ell (K_{r,n}-1)}{2n} \right) \nonumber \\
&\leq \mathrm{exp} \left(\frac{1-\tilde{\mu}}{\tilde{\mu}} \ell \left(\frac{1}{2} \right)^{K_{r,n}-1} -\left( 1-\tilde{\mu} \right) \frac{\ell (K_{r,n}-1)}{4} \right) \label{eq:StepbyStep1} \\
&= \mathrm{exp} \left(-\left( 1-\tilde{\mu} \right) \ell \left(\frac{(K_{r,n}-1)}{4} - \frac{\left(0.5\right)^{K_{r,n}-1}}{\tilde{\mu}} \right) \right) \nonumber \\
&\leq \mathrm{exp} \left( -2\left( 1-\tilde{\mu} \right) \left(\frac{(K_{r,n}-1)}{4}  - \frac{(0.5)^{K_{r,n}-1}}{\tilde{\mu}} \right) \right) \label{eq:StepbyStep3}
\end{align}
where (\ref{eq:StepbyStep1}) follows from the facts that $n-\ell \geq n/2$ and $\ell/n \leq 0.5$ on the specified range for $\ell$, and  (\ref{eq:StepbyStep3}) follows 
for all $K_{r,n}$ sufficiently large such that
$
K_{r,n} \geq \ceil*{4 \left( \frac{(0.5)^{K_{r,n}-1}}{\tilde{\mu}}\right)+1}
$ upon noting that $\ell \geq 2$.

If, on the other hand, it holds that  $\frac{\ell(K_{r,n}-1)}{n} > 1$, we see that
$1-e^{-\ell \left( \frac{K_{r,n}-1}{n}\right)} \geq 1-  e^{-1}$. 
Reporting this into  (\ref{eq:defn_A_nr}) 
and using $\ell \leq n/2$, we get
\begin{align}\nonumber
A_{n,\ell} &\leq \mathrm{exp} \left(\frac{1-\tilde{\mu}}{\tilde{\mu}} \ell \left(\frac{\ell}{n} \right)^{K_{r,n}-1} -\left( 1-\tilde{\mu} \right) (n-\ell) \left(1-  e^{-1} \right)  \right) \nonumber \\
&\leq \mathrm{exp} \left(\frac{1-\tilde{\mu}}{\tilde{\mu}} \left( \frac{n}{2} \right) (0.5)^{K_{r,n}-1} -\left( 1-\tilde{\mu} \right) \frac{n}{2} \left(1-  e^{-1} \right)  \right) \nonumber \\
&= \mathrm{exp} \left(-\left( 1-\tilde{\mu} \right) \frac{n}{2} \left( 1-e^{-1}  - \frac{(0.5)^{K_{r,n}-1}}{\tilde{\mu}} \right)  \right).  \label{eq:StepbyStep4}
\end{align}

Combining (\ref{eq:StepbyStep3}) and (\ref{eq:StepbyStep4})
 we see that $A_{n, \ell} \leq \Psi(n, \pmb{\mu}, \pmb{K}_n) $ for all $n$ sufficiently large and all $\ell=2, \ldots, \lfloor \frac{n}{2} \rfloor$, where $\Psi(n, \pmb{\mu}, \pmb{K}_n)$ is given by (\ref{eq:defining_Psi}).
 
Observing that the bound derived on $A_{n,\ell}$ is independent on $\ell$, 
we get from (\ref{eq:defining_Psi}) and (\ref{eq:last_bound_osy}) 
\begin{equation}\nonumber
\sum_{\ell=2}^{\lfloor \frac{n}{2} \rfloor} \binom{n}{\ell} \mathbb{P}\left[ B_{n,\ell}(\pmb{\mu},\pmb{K}_n ) \right] \leq  \Psi(n, \pmb{\mu}, \pmb{K}_n)  \sum_{\ell=2}^{\infty}  \tilde{\mu}^\ell = \frac{\tilde{\mu}^2}{1-\tilde{\mu}} \Psi(n, \pmb{\mu}, \pmb{K}_n) 
\end{equation}

Letting $n$ go to infinity, it is now easy to see that
\[
\lim_{n \to \infty} \Psi(n, \pmb{\mu}, \pmb{K}_n) = 0, \quad 2 \leq \ell \leq n/2
\]
under the enforced assumption that $\lim_{n \to \infty} K_{r,n}= \infty$. Hence, the conclusion
\begin{equation}\nonumber
\lim_{n \to \infty} \sum_{\ell=2}^{\lfloor \frac{n}{2} \rfloor} \binom{n}{\ell} \mathbb{P}\left[ B_{n,\ell}(\pmb{\mu},\pmb{K}_n ) \right] = 0
\end{equation}
immediately follows since $0 < \tilde{\mu} <1$. This establishes Theorem~\ref{thm:OneLaw+Connectivity}.

\section{Conclusion}
In this paper, we have proposed inhomogeneous random K-out graphs $\mathbb{H}\left(n;\pmb{\mu},\pmb{K}_n\right)$ where nodes are arranged into $r$ disjoint classes and the number of selections made by a node is dependent on its class. In particular, we consider the case where each node is classified as class-$i$ with probability $\mu_i>0$ for $i=1,\ldots,r$. A class-$i$ node selects $K_{i,n}$ other nodes uniformly at random to be paired to. Two nodes are deemed adjacent if at least one selects the other. Without loss of generality, we assumed that $K_{1,n} \leq K_{2,n} \leq \ldots \leq K_{r,n}$. 

Earlier results on homogeneous random K-out graphs (where all nodes select $K$ other nodes) suggest that the graph is connected whp if $K\geq2$. Hence, $\mathbb{H}\left(n;\pmb{\mu},\pmb{K}_n\right)$ is trivially connected whenever $K_{1,n} \geq 2$. We investigated the connectivity of $\mathbb{H}\left(n;\pmb{\mu},\pmb{K}_n\right)$ in the particular case when $K_{1,n}=1$. Our results revealed that when $K_{1,n}=1$, $\mathbb{H}\left(n;\pmb{\mu},\pmb{K}_n\right)$ is connected with high probability if and only if $K_{r,n} = \omega(1)$. Any bounded choice of $K_{r,n}$ is shown to yield a positive probability of $\mathbb{H}\left(n;\pmb{\mu},\pmb{K}_n\right)$ being {\em not} connected, and an explicit lower bound on this probability is provided.

\bibliographystyle{IEEEtran}
\bibliography{IEEEabrv,references}

\end{document}